\documentclass{tran-l}

 \newtheorem{thm}{Theorem}[section]
 
 \newtheorem{lem}[thm]{Lemma}
 
 \theoremstyle{definition}
 \newtheorem{defn}[thm]{Definition}
 \theoremstyle{remark}
 \newtheorem{rem}[thm]{Remark}
 \numberwithin{equation}{section}
 \newcommand{\Real}{\mathbb{R}}

\begin{document}

\title[Arithmetic Operations and Generalized Derivatives]
 {The Natural Chain of Binary Arithmetic Operations and Generalized Derivatives}

\author{ Michael L. Carroll }

\address{BEI Systron Donner Inertial Division, Concord, CA 94518
 }

\email{mcarroll@systron.com}



\date{TBD and, in revised form, TBD}

\dedicatory{}

%%% ----------------------------------------------------------------------

\begin{abstract}
This paper presents a graded hierarchy or chain of binary
operations on the reals and the complex numbers. The operations
are related distributively in the sense that any one of them
distributes over the next lower operation in the chain. For one
particular operation we explore specific properties and derive
results, including several useful formulas and identities. Next,
the operation is extended to the complex numbers and a new kind of
derivative is defined based on this binary operation.  Some basic
formulae analogous to standard classical ones are proven for this
new derivative.  Finally, the derivative is generalized to the
point that the new derivative, the classical one, and countably
many others are seen to be special cases.
\end{abstract}

%%% ----------------------------------------------------------------------
\maketitle
%%% ----------------------------------------------------------------------

\section*{Introduction}

Nearly every culture in the world has discovered and made use of
the natural numbers and their extensions the integers, the
rationals, and beyond.  The natural operation of counting gives
rise quite naturally to the operations of addition and
multiplication. Moreover, these operations are intimately related.
Indeed, multiplication is usually defined in terms of addition and
is considered universally as a kind of shorthand for addition.
Perhaps because of this, multiplication and addition have a
natural relationship expressed by the distributive law. This law
binds these two operations together and endows the set of integers
with the algebraic structure of a ring.

The logarithm and exponential maps are compatible with the
intimate relationship between multiplication and addition.  The
logarithm maps products into sums in such a way that the
multiplicative structure of the positive real numbers and the
additive structure of all the reals are in essence mirror images
of each other; as groups they are isomorphic.

Are addition and multiplication the only natural operations on the
natural numbers whose relationship is expressed by a distributive
law and a mapping like the logarithm?  Are there other natural
arithmetic operations hidden in the integers or reals that are
related in a similar way?

For a long time algebraists have studied other ring and field
structures having operations that are related to each other as
addition is to multiplication. In fact, the linear, homogenous
nature of ring endomorphisms suffices to endow any Abelian group
with a natural ring structure. (See \cite{Hun74}, pg. 116.)

One group of algebraists has paid considerable attention to
max-plus algebras and semirings (See \cite{Gau97}).  These
structures usually employ the operation of addition and the
maximum operator: $\max\{x,y\}$.  In this case, ordinary addition
distributes over $\max$. Results in this area have contributed
greatly to the solution of problems in optimization and control.

The purpose of the present paper is to show that addition and
multiplication are two adjacent operations in an infinite chain or
graded hierarchy of binary operations on the reals and complex
numbers. The chain can be extended "upwards" from multiplication
and "downwards" from addition.  Each operation is related to its
next lower neighbor in the chain by being distributive over its
neighbor.  The logarithm provides the means for descending down
the chain and the exponential function the means for ascending up
the chain.

After we exhibit some general algebraic properties for these
operations, we will proceed to investigate one of these new
operations -- the join operation -- in considerably more detail.
We shall discover that when extended to the complex numbers the
join / additive structure is potentially as rich as the ordinary
additive / multiplicative structure of the complex numbers.

%%% ----------------------------------------------------------------------

\section{The Natural Chain of Binary Arithmetic Operations}

Throughout this paper $\Real$ will denote the field of real
numbers, $\Real_0^+$ the multiplicative monoid of positive reals
with 0 adjoined, and $\Real_\infty$ the extended real number
system, i.e., the real numbers with the two extreme points
$\pm\infty$ adjoined and the usual operational conventions.
$\Real_\infty$ is not a field, however, since expressions like
$\infty-\infty$ and $\infty/\infty$ are not defined.  Further, we
will extend the definition of the natural logarithm to include 0
in its domain: $\log(0):=-\infty$. Similarly, we define
$e^{-\infty}:=0$, and $e^{+\infty}:=\infty$.

With these conventions we now define a countable family of binary
operations on $\Real_\infty$ as follows:

\begin{defn}
Let $x,y\in\Real_\infty$.  For each $n\in\mathbb{Z}$, we shall
recursively define a binary operation $\oplus_n$ on $\Real_\infty$
as follows.  For $n=0$, define
\begin{equation}
x\oplus_0 y:=x+y
\end{equation}
For each $n\leq 0$, define
\begin{equation}
\label{eqno2} x\oplus_{n-1} y:=\ln(e^{x}\oplus_n e^{y})
\end{equation}
Finally, for $n\geq 0$, define
\begin{equation}
\label{eqno3} x\oplus_{n+1} y:=\exp[\ln(x)\oplus_n \ln(y)]
\end{equation}
provided $x\geq 0$ and $y\geq 0$. We shall call $\oplus_n$ an
\textbf{\textit{arithmetic operation}}.
\end{defn}

\begin{rem}
We note that for ordinary multiplication ($\times$) we have
$\times=\oplus_1$. Also, we introduce the notation $\vee$ for the
special case of $n=-1$: i.e., in place of $\oplus_{-1}$ we shall
use $\vee$. Two operations are said to be
\textbf{\textit{adjacent}} if their associated integers are
consecutive. Adjacent operations are, as we shall see below, very
closely related.

Many authors, when considering only two adjacent operations, tend
to use $\otimes$ for the "higher" operation and $\oplus$ for the
lower one.  However, since in this paper we will often consider
three operations together, we will either retain the general
notation of $\oplus_n$ or use $\vee$ in conjunction with $+$ and
$\times$ (the latter symbol being usually omitted), depending on
the context. The choice of the $\vee$ notation is motivated by its
use in lattice theory for an addition-like operation. It is read
as "join," "cup," or "vee" (although "union" might also be
appropriate).  The join terminology, which I tend to favor, is
taken from lattice theory \cite{Gra68}. In our new algebraic
system $\vee$~ will play an additive role while $+$ will play a
multiplicative role.

Eq. ~(\ref{eqno2}) of this definition makes sense because $e^x\geq
0$ for all $x\in\Real_\infty$. The restriction on Eq.
~(\ref{eqno3}) prevents us from applying the logarithm to negative
numbers, a restriction we shall remove when we extend the domain
of definition to the complex numbers.
\end{rem}

\begin{rem}
Eq. ~(\ref{eqno2}) appears to be a special case of an operation
defined by Maslov et. al. in \cite{Mas95} and referenced also by
Gaubert in \cite{Gau97}.  Maslov defines: $x\oplus_h y:=
h\ln(e^{x/h}+e^{y/h})$.  By letting $h=1$ in Maslov's definition
we obtain what I have called $\oplus_{-1}$.  The semirings
resulting from use of $\oplus_{h}$ and addition according to
Maslov's definition are sometimes referred to as log-plus
semirings (or semi-fields when extended with inverses, see
\cite{McE00}). The Maslov theory is primarily concerned with a
single operation in relation to addition, parameterized by the
continuous parameter $h$ with convergence to the max operation in
the limit as $h\to 0$. The theory presented here, however, is
concerned with distinct operations for each $n$, and with the
special cases of addition when $n=0$ (by convention) and
multiplication when $n=1$.
\end{rem}

Using the notation $f^{(n)}:=\mathop{\underbrace{f\circ f \circ
\cdots \circ f }_{n-times}}$, where $\circ$ represents function
composition, it is easily verified that for $n>0$
\begin{equation}
\label{up} x\oplus_n y= \exp^{(n)}[\ln^{(n)}(x)+\ln^{(n)}(y)]
\end{equation}
for non-negative $x,y\in\Real_\infty$ and
\begin{equation}
\label{down} x\oplus_{-n} y=\ln^{(n)}[\exp^{(n)}(x)+\exp^{(n)}(y)]
\end{equation}
for all $x,y\in\Real_\infty$.

%%% ----------------------------------------------------------------------

\section{Algebraic Properties of Arithmetic Operations}

We now explore a few basic algebraic properties of adjacent
operations.  We will focus primarily on traversing the chain
downward:

\begin{thm}  On $\Real_\infty$, for all $n\in\mathbb{Z}$,

(i)   $\oplus_n$ is associative.

(ii)  $\oplus_n$ is commutative.

(iii) $\oplus_n$ distributes over $\oplus_{n-1}$.
\end{thm}
\begin{proof}
First, we note that Eq. (\ref{up}) implies that for $x,y\geq 0$,
$\ln(x\oplus_n y) = \ln(x)\oplus_{n-1} \ln(y)$.  For,
\begin{align*}
\ln(x\oplus_n y) & = \ln\circ \exp^{(n)}[\ln^{(n)}(x)+\ln^{(n)}(y)] \\
& = \exp^{(n-1)}[\ln^{(n)}(x)+\ln^{(n)}(y)] \\
& = \exp^{(n-1)}[\ln^{(n-1)}\circ\ln(x)+\ln^{(n-1)}\circ\ln(y)] \\
& =  \ln(x)\oplus_{n-1} \ln(y),
\end{align*}
(i) To prove associativity for $n<0$, assume inductively that
$\oplus_n$ is associative and consider:
\begin{align*}
  (x\oplus_{n-1} y) \oplus_{n-1} z & = \ln{ (e^{\ln{ (e^{x} \oplus_{n} e^{y}) }} \oplus_n e^{z})} \\
  & = \ln{(e^{x}\oplus_n e^{y} \oplus_n e^{z})} \\
  & = \ln {(e^{x} \oplus_n e^{\ln{(e^{y} \oplus_n e^{z})}} )} \\
  & = x\oplus_{n-1} (y \oplus_{n-1} z).
\end{align*}
Similarly, for $n>0$, and subject to the restriction $x,y\geq 0$,
we have
\begin{align*}
  (x\oplus_{n+1} y) \oplus_{n+1} z & = \exp[ \ln(x \oplus_{n+1} y) \oplus_{n} \ln{z}] \\
  & = \exp[\ln{x}\oplus_{n} \ln{y} \oplus_{n} \ln{z}] \\
  & = \exp[\ln{x} \oplus_{n} \ln(y \oplus_{n+1} z)] \\
  & = x\oplus_{n+1} (y \oplus_{n+1} z).
\end{align*}
(ii) This is an obvious consequence of the commutativity of
addition. \\
(iii) To prove the distributive property, consider
\begin{align*}
  x \oplus_n (y \oplus_{n-1} z) & =
  \exp^{(n)}[\ln^{(n)}x+\ln^{(n)}(y\oplus_{n-1} z)] \\
  & =    \exp^{(n)}[\ln\circ\ln^{(n-1)}x+\ln\circ\ln^{(n-1)}(y\oplus_{n-1} z)] \\
  & =    \exp^{(n)}[\ln(\ln^{(n-1)}x\times\ln^{(n-1)}(y\oplus_{n-1} z))] \\
  & =    \exp^{(n)}[\ln(\ln^{(n-1)}x\times [\ln^{(n-1)}y+\ln^{(n-1)}z])] \\
  & =    \exp^{(n)}[\ln(\ln^{(n-1)}x\times \ln^{(n-1)}y+\ln^{(n-1)}x\times \ln^{(n-1)}z)] \\
  & =    \exp^{(n)}\circ\ln[\ln^{(n-1)}x\times \ln^{(n-1)}y+\ln^{(n-1)}x\times \ln^{(n-1)}z] \\
  & =    \exp^{(n-1)}[\ln^{(n-1)}x\times \ln^{(n-1)}y+\ln^{(n-1)}x\times \ln^{(n-1)}z] \\
  & = (x \oplus_n y) \oplus_{n-1} (x \oplus_n z).
\end{align*}
\end{proof}
Note that $-\infty$ is an identity element for $\vee$, because
 \[ x\vee -\infty = \ln{(e^{x}+e^{-\infty})} = \ln{(e^{x}+0)} =
\ln{(e^{x})} = x. \] It is also two-sided owing to the commutative
property.  Going up the chain we discover a very interesting
sequence of identity elements and inverse elements. Table
\ref{tbl1} lists some of the identities and inverses for these low
order operations.
\begin{table}[h]
\label{tbl1}
\begin{tabular}{|r|c|c|c|}
\hline
\textbf{Value of n} & \textbf{Common Notation} & \textbf{Identity} & \textbf{Inverse} \\
\hline
$n=3$ & N/A & $e^e$ & $\exp[\exp(1/\ln\ln{x})]$ \\
\hline
$n=2$ & N/A & $e$ & $\exp(1/\ln{x})$ \\
\hline
$n=1$ & $\times$ & $1$ & $1/x$ \\
\hline
$n=0$ & $+$ & $0$ & $-x$ \\
\hline
$n=-1$ & $\vee$ & $-\infty$ & TBD \\
\hline
\end{tabular}
\caption{Operations, Identities and Inverses}
\end{table}

The TBD (=To Be Determined) will be cleared up in the complex
case.  The operation $\oplus_2$ captures the essence of
exponentation. For, we note that $e^{n}\oplus_2 x=x^{n}$, and
$e^{m}\oplus_2 e^{n}=e^{mn}$, for any $x\in\Real_\infty$ and
$m,n\in\mathbb{Z}$. Incidently, the commutativity of $\oplus_2$
implies the following beautiful and well-known fact:
\begin{equation}
x^{\ln y}=y^{\ln x}.
\end{equation}

%%% ----------------------------------------------------------------------

\section{Algebraic Properties of the Join Operation}

In this section we will begin our detailed study of the join
operation $\vee$, defined for all $x,y\in\Real_\infty$. If you
think of multiplication as being "fast addition", then the $\vee$
operation can be thought of as "slow addition." Because it is the
image of a sum under the logarithm mapping, it could also be
called "logarithmic addition."

For positive $x\in\mathbb{N}$ we use the notation
$nx=\mathop{\underbrace{x + \cdots + x }_{n-times}}$.  Similarly,
we introduce the following notation for $\vee: n_\vee
x:=\mathop{\underbrace{x \vee \cdots \vee x }_{n-times}}$.  We
also agree that $1_\vee x=x$.  With this notation we note that

\begin{lem}
\begin{equation}\label{eq3-7}
  n_\vee x= \ln(n)+x
\end{equation}

\end{lem}
\begin{proof}
We compute:

$n_\vee x:=\mathop{\underbrace{x \vee \cdots \vee
x}_{n-times}}=\ln(e^{x}+\cdots
+e^{x})=\ln(ne^{x})=\ln(n)+\ln(e^{x})=\ln(n)+x$.
\end{proof}

An interesting application of this result is
$\mathop{\underbrace{0 \vee \cdots \vee 0 }_{n-times}}=\ln(n)$.
Also note that $\ln(nx)=n_\vee\ln(x)$.  This is a consequence of
the homomorphic property of the logarithm: $\ln(x+y)=\ln(x)\vee
\ln(y)$. This homomorphism is in fact an isomorphism with the
$\exp$ function as the inverse map.  Thus we have $e^{x\vee
y}=e^x+e^y$. It is also easy to verify the following two equations
for integers $m$ and $n$ :
\begin{equation}\label{eq3-8}
  (mn)_\vee x=m_\vee (n_\vee x)
\end{equation}
\begin{equation}\label{eq3-9}
  m_\vee x +n_\vee y= \ln(mn)+x+y
\end{equation}
Some useful recurrence relations are the following:
\begin{equation}\label{eq3-10}
  m\vee n =[(m-1)\vee (n-1)]+1
\end{equation}
\begin{equation}\label{eq3-11}
  (m\vee n)+1= (m+1)\vee (n+1)
\end{equation}

One of the interesting properties of the Binomial Theorem is that
it relates operations that are not adjacent but are once removed,
so to speak.  The ordinary result from algebra, for instance,
indicates what happens when you apply integer exponentiation to a
sum.  We now wish to prove for our new system the analogue of this
theorem by considering what happens when we apply integer
multiplication to a join. For $n=2$ we compute:
\begin{equation}\label{binom2}
  2(x\vee y)=2x\vee 2_\vee (x+y)\vee 2y
\end{equation}
Here we have used the following:
\begin{align*}
2(x\vee y)&=(x\vee y)+(x\vee y) \\
&=[x+(x\vee y)]\vee [y+(x\vee y)] \\
&=(x+x) \vee (x+y) \vee (y+x) \vee (y+y) \\
&=2x \vee 2_\vee (x+y) \vee 2y \\
&=1_\vee 2x \vee 2_\vee (x_y) \vee 1_\vee 2y,
\end{align*}
where we have used the convention that $1_\vee x = x$. A couple of
useful formulas are the following:
\begin{equation}\label{eq3-13}
  n_\vee (x+y) + z=n_\vee (x+y+z)
\end{equation}
\begin{equation}\label{eq3-14}
  x + y + (x\vee y) = (2x+y)\vee (x+2y)
\end{equation}
 We also introduce the notation $\bigvee_{k=0}^n
x_{k}= x_0\vee \cdots \vee x_{n}$.
\begin{thm}(Binomial Theorem)
\begin{equation}\label{binomn}
  n(x\vee y)=\bigvee_{k=0}^{n} {n \choose k}_\vee [(n-k)x+ky]
\end{equation}
\end{thm}
\begin{proof}
Using the fact that $\exp(n_\vee x)=n\exp(x)$, and making the
substitutions $u=\exp(x)$ and $v=\exp(y)$, we see that
\begin{align*}
\exp\{\bigvee_{k=0}^{n} {n \choose k}_\vee [(n-k)x+ky]\} &= \sum_{k=0}^{n} {n \choose k} \exp[(n-k)x+ky] \\
&= \sum_{k=0}^{n} {n \choose k} \exp[(n-k)x]\exp[ky] \\
&= \sum_{k=0}^{n} {n \choose k} u^{n-k}v^{k} \\
&= (u+v)^n.
\end{align*}
Taking the logarithm of both sides yields the desired formula.
\end{proof}

Although we will not explore polynomials to any extent in this
paper, we mention in passing that a general polynomial using
addition and join can be written thus:
\begin{equation}
 p(x)=(a_n+nx)\vee (a_{n-1}+(n-1)x)\vee \cdots \vee (a_1+x)\vee a_0
\end{equation}
By this time the reader should be comfortable enough with the
notation to allow us to drop some of the parentheses, bearing in
mind that $+$ binds more closely than $\vee$.  Therefore, we can
declutter the previous equation as:
\begin{equation}
 p(x)=a_n+nx\vee a_{n-1}+(n-1)x\vee \cdots \vee a_1+x\vee a_0
\end{equation}
It should be quite apparent by now that many of the ordinary
arithmetic expressions and statements involving exponentiation,
multiplication and addition can be readily transformed into ones
involving multiplication, addition, and join, respectively.  As we
saw in the proof of the Binomial Theorem, many analogues of
standard theorems can be successfully guessed by taking
logarithms, and their proofs can be transferred back into the
domain of ordinary algebra by applying the exponential function.
It is this transfer principle that gives rise to so many
interesting results using join and addition. It should also be
borne in mind that many of the results we are demonstrating can be
obtained easily for arbitrary pairs of adjacent operations. The
algebra becomes a bit more complex, however, and we will leave
exploration of more general cases to subsequent papers.

For now we will leave algebra behind and begin to explore the
analytic properties of the join operation.

\section{Analytic Properties of the Join Operation}
The join operation provides a smooth approximation of the $\max$
operation.  The following facts should aid in making this clear
and helping us to understand the join's qualitative behavior:
$\lim_{x\to -\infty} (x\vee y)=\lim_{x\to -\infty} \ln(e^x + e^y)
= y$. This implies that $x\vee y \approx y$ whenever $x$ is
significantly smaller than $y$.  On the other hand, $x\vee y
\approx x$, whenever $x$ is considerably larger than $y$.  Thus
the value of the join depends on which exponential term dominates.

The join operation, as a composite of continuous functions is
itself continuous at every point in $\Real_\infty$.  Moreover, it
is clearly differentiable in both arguments.  We will now explore
how it behaves under ordinary and partial differentiation.
\subsection{Differential Identities and Formulas}

We wish to determine how $\vee$ behaves under differentiation.  If
$u=u(x)$ and $v=v(x)$ are differentiable real functions of a
single real variable, then the elementary formulas of the
differential calculus tell us that
\begin{equation}
  \frac{d}{dx}uv=u\frac{dv}{dx}+\frac{du}{dx}v
\end{equation}
\begin{equation}
  \frac{d}{dx}(u+v)=\frac{du}{dx}+\frac{dv}{dx}
\end{equation}

An operator such as $d/dx$ on a semiring R with these properties
is called a derivation on R.  We can see easily enough that $d/dx$
is \textit{not} a derivation on the semiring $\langle\Real_\infty
; +,0 ; \vee , -\infty\rangle$.  In fact, we can easily compute
the general formula:
\begin{equation}\label{dervofjoin}
  \frac{d}{dx}(u\vee v)= \frac{e^{u}\displaystyle \frac{du}{dx}+e^{v}\displaystyle \frac{dv}{dx}}{e^u + e^v}
\end{equation}
In particular, we see that
\begin{equation}\label{dq4-18}
  \frac{d}{dx}(u\vee u)= \frac{e^{u}\displaystyle\frac{du}{dx}+e^{u}\displaystyle\frac{du}{dx}}{e^u +
  e^u}=\frac{2e^u}{2e^u}\frac{du}{dx}=\frac{du}{dx}.
\end{equation}
If $u(x)=x$, then $\displaystyle \frac{d}{dx}(x\vee x)=1$, from
which it follows that the function $x\vee x=x+C$, which is
consistent with what we already know, i.e., that $x\vee
x=x+\ln(2)$, or $C=\ln(2)$.

Also note that $\displaystyle \frac{d}{dx}(x\vee -x)=\tanh(x)$,
from which we derive $x\vee -x=\ln \cosh(x)+C$, to within a
constant. Note that $\ln \cosh(0)=0$ and $\ln(2)=0\vee 0$, which
implies that $x\vee -x=\ln \cosh(x)+\ln(2)=\ln(2 \cosh(x))$.  This
is consistent, of course, with our original definition of the
join. We'll obtain even more interesting results along these lines
when we extend the join to the complex numbers.  Also, Eq.
~(\ref{dervofjoin}) will turn out to be of considerable utility in
demonstrating many properties of a new derivation we will define.

We close this section with a few more formulae, the details of the
proofs being simple and therefore omitted:
\begin{thm}
\begin{equation}\label{eq4-5}
  \left [\frac{\partial}{\partial x}+\frac{\partial}{\partial y}\right ](x\vee y)=1
\end{equation}
\begin{equation}\label{eq4-6}
  \nabla^2(x\vee y)=1 - \left [\left(\frac{\partial}{\partial x}(x\vee y)\right )^2+\left (\frac{\partial}{\partial y}(x\vee y)\right)^2\right ]
\end{equation}
\begin{equation}\label{eq4-7}
  \frac{\partial}{\partial x}\frac{\partial}{\partial y}(x\vee y)=  \frac{\partial}{\partial y}\frac{\partial}{\partial x}(x\vee
  y)=\frac{-e^xe^y}{e^x + e^y}=-\exp\left [ x+y - (x\vee y) \right ]
\end{equation}
\begin{equation}\label{eq4-8}
  \left [ \frac{\partial}{\partial x} \vee \frac{\partial}{\partial y} \right ](xy)=x\vee y
\end{equation}
\begin{equation}\label{eq4-9}
  \left [ \frac{\partial}{\partial x} \vee \frac{\partial}{\partial y} \right
  ](x+y)=1\vee 1=1+\ln2
\end{equation}
\begin{equation}\label{eq4-10}
  \left [ \frac{\partial}{\partial x} \vee \frac{\partial}{\partial y} \right ](x\vee
  y)=\frac{e^x}{e^x+e^y}\vee \frac{e^y}{e^x+e^y}
\end{equation}
\end{thm}

Although these formulae are interesting in their own right and
potentially useful in simplifying calculations and solving
differential equations, what we really need is a ring (or better
yet, a field) and a new derivation that parallels the classical
differentiation.  To accomplish this we must now turn to the
complex numbers.
\section{Extension of the Join to the Complex Plane}

Extending the join operation to the complex plane is
straightforward and has several benefits, the chief of which being
that we obtain a full and rich algebraic structure that allows us
to construct a new differential and integral calculus.

Before heading down that path, however, we must first take a
little care in dealing with the point at infinity.  In extended
real analysis, we distinguish between $+\infty$ and $-\infty$.
They are two different points added to the real line. In the
complex plane, however, thanks to the Riemann Sphere, we usually
consider only one point at infinity. This is inconvenient when we
wish to have $e^{-\infty}=0$. Although this statement makes
intuitive sense in the extended real system, it is not so
intuitive in the complex case.  For, in the complex case all
points at infinity tend to coalesce into a single point.

There are several ways to get around this issue.  First, we can
distinguish in terms of how the point at infinity is approached.
Left and right hand limits are in some sense different.  For
example, if we treat an approach to $\infty$ with the real part
becoming more and more negative, we might call this approaching
$\infty$ from the right.  Similarly, approaching $\infty$ from the
left would mean that the real part is becoming more and more
positive.  This would be consistent with $\lim_{z\to
-\infty}e^z=0$ and $\lim_{z\to \infty}e^z=\infty$.  Obviously,
since these two values are not the same, if we treat $+\infty$ and
$-\infty$ as the same point at infinity, then the exponential
function has a discontinuity at that point.  However, if we keep
these points distinct, we need not worry about the limits being
different.

Another approach would be to adjoin four new points to
$\mathbb{C}$:  The real infinities $+\infty$ and $-\infty$ and the
pure imaginary infinities $+i\infty$ and $-i\infty$.  The
algebraic operations using these symbols will make sense as long
as we don't try to cancel infinities in undefined ways.  This
approach destroys, however, the field structure of the complex
numbers.

For our purposes it is not necessary to proceed in this manner. We
can retain a field structure on $\mathbb{C}$ by adjoining a single
new infinite quantity: $-\infty$.  Thus $\langle
\mathbb{C}\cup{-\infty};\vee, -\infty;+,0\rangle$ becomes a field
with the conventions that $-\infty+z=-\infty$ and $-\infty \vee
z=z$ for all $z\in\mathbb{C}$.  Addition is the "multiplicative"
operation of the field and join is the "additive" operation.

To complete the field structure we need to show that the join
operation in the complex plan affords an "additive" inverse.
However, before we do this, we need to discuss the multifunction
nature of the complex logarithm.

Since the complex logarithm is a multifunction that takes on
infinitely many values for each point of the complex plane, we
need to take care when asserting equality of expressions involving
logarithms. Because the principal branch of the logarithm Log($z$)
is not a homomorphism (i.e., it does not always map products into
sums), we will avoid limiting ourselves to this branch. Therefore,
we shall use the general complex logarithm
$\log(z)=\ln|z|+i\arg(z)$ and be deliberately vague about which
branch of the logarithm we restrict ourselves to in any particular
situation.  In many cases our equations will, strictly speaking,
be equivalences modulo $2\pi i$.  For further discussion of
various branches of the logarithm, see Ahlfors \cite{Ahl66}.

We can now define our complex join as:
\begin{equation}\label{complexjoin}
  z_1\vee z_2 := \log\left [\exp(z_1) + \exp(z_2) \right ]
\end{equation}

By the convention we have adopted regarding the absorptive nature
of $-\infty$, we can see that $\log(z) \to -\infty$ continuously
as $z \to 0$.

The inverse of $z$ with respect to the join operation will
sometimes be denoted by $\tilde{z}$.  However, we can readily
compute that
\begin{equation}
\tilde{z}=z+i\pi
\end{equation}
For,
\begin{align*}
  z\vee (z+i\pi) &= \log(e^z + e^{z+i\pi}) \\
  &= \log(e^z+e^{z}e^{i\pi}) \\
  &= \log(e^z - e^z) \\
  &= \log(0) \\
  &= -\infty.
\end{align*}
We also notice that $\log(-z)=\log(z) +i\pi$, for
\begin{align*}
  \log(-z)\vee \log(z) &= \log(e^{\log(-z)}+e^{\log(z)}) \\
   &= \log(-z+z) \\
   &= \log(0) \\
   &= -\infty.
\end{align*}
With the definition of the inverse in hand, we are now in a
position to define a new kind of derivative:
\begin{defn}
Let $G\subseteq\mathbb{C}$ be an open subset of $\mathbb{C}$,
$f:G\to \mathbb{C}\cup\{-\infty\}$, and suppose that for some
$z\in\mathbb{C}\cup\{-\infty \}$ the limit
\begin{equation}\label{newderiv}
  \lim_{h\to -\infty} \left \{ \left [ f(z\vee h) \vee (f(z)+i\pi) \right ] - h \right \}
\end{equation}
exists and is finite.  Then we shall say that $f$ is
$\vee$-\textit{\textbf{differentiable}} at $z$ and denote the
limit by $D_\vee f(z)$, called the
$\vee$-\textbf{\textit{derivative}} of $f$ at $z$.
\end{defn}
\begin{rem}
Whenever there is little danger of confusion and the meaning is
clear from the context, we shall simply write $Df(z)$ instead of
$D_\vee f(z)$.  Later, when we once again consider more general
operations on $\mathbb{C}$ and define more generalized
derivatives, we will use the notation $D_{-1}$ instead of either
$D$ or $D_\vee$.
\end{rem}
\begin{rem}
This definition results from a direct translation of the ordinary
difference quotient $\displaystyle \frac{f(z+h)-f(z)}{h}$ using
the previously mentioned transfer principle that sends sums into
joins, subtractions into joins of "additive" inverses, and
quotients into subtraction.  Thus this definition is the natural
analogue of the ordinary derivative.
\end{rem}

Some basic properties of $\vee$-differentiation that are near
obvious analogues of their classical counterparts are the
following:

\begin{thm}
\label{bigthm} Let $n$ be a non-negative integer, $a$ a constant
complex number, $z$ an indeterminate complex variable, and both
$f$ and $g$ complex functions analytic and $\vee$-differentiable
on an open set in $\mathbb{C}$.  Then \\
\flushleft
\begin{itemize}
\item[(i)] $D(a)=-\infty$
\item[(ii)] $D(z)=0$
\item[(iii)] $D(a+z)=a$
\item[(iv)] $D(n_\vee f)=n_\vee D(f)$
\item[(v)] $D(f\vee g)=D(f)\vee D(g)$
\item[(vi)] $D(f+g)=[f+D(g)]\vee [D(f)+g]$
\item[(vii)] $D(nz) = n_\vee (n-1)z=(n-1)z+\ln(n)$
\item[(viii))] $D(z^{n})= \log \left[ \displaystyle \frac{nz^{n-1}e^{z^{n}}}{e^z} \right] =z^n -z+(n-1)\log(z)+\ln(n)$
\item[(ix)] $D^{n}(nz)= \ln 1+\cdots + \ln n$
\end{itemize}
\end{thm}
\begin{proof}
Although these proofs are elementary, they serve to illustrate
some important techniques that are peculiar to reasoning with
joins. Later, we shall see that each of these formulae can be
simply derived from a single principle expressing the relationship
between $\vee$-differentiation and ordinary differentiation.
\flushleft
\begin{itemize}
\item[(i)]  Observe that
\begin{align*}
D(a)&= \lim_{h \rightarrow -\infty} \{[a\vee (a+i\pi ) ] -h\}\\
&=\lim_{h \rightarrow -\infty} \{-\infty -h\}\\
&=\lim_{h \rightarrow -\infty} \{-\infty\}\\
&=-\infty.
\end{align*}
\item[(ii)] Again, we compute
\begin{align*}
  D(z) &=\lim_{h \rightarrow -\infty} \{[z\vee h\vee (z+i\pi)]-h\} \\
  &= \lim_{h \rightarrow -\infty} \{(z-h)\vee (h-h)\vee (z+i\pi -h)\}\\
  &= \lim_{h \rightarrow -\infty} \{(z-h)\vee 0 \vee (z+i\pi -h)\}\\
  &= \lim_{h \rightarrow -\infty} \{(z-h)\vee (z-h+i\pi)\vee 0\}\\
  &= \lim_{h \rightarrow -\infty} \{-\infty\vee 0\}\\
  &= \lim_{h \rightarrow -\infty} \{0\}\\
  &=0.
\end{align*}
\item[(iii)] Similarly
\begin{align*}
  D(a+z) &=\lim_{h \rightarrow -\infty} \{[a+(z\vee h)\vee (a+z+i\pi)]-h\} \\
  &= \lim_{h \rightarrow -\infty} \{[(a+z)\vee (a+h)\vee (a+z+i\pi)] -h\}\\
  &= \lim_{h \rightarrow -\infty} \{(a+z-h)\vee (a+h-h) \vee (a+z+i\pi -h)\}\\
  &= \lim_{h \rightarrow -\infty} \{(a+z-h)\vee (a+z-h+i\pi)\vee a\}\\
  &= \lim_{h \rightarrow -\infty} \{-\infty\vee a\}\\
  &= \lim_{h \rightarrow -\infty} \{a\}\\
  &=a.
\end{align*}
\item[(iv)] First note that
\begin{align*}
   (f\vee g)(z\vee h) \vee ([(f\vee g)(z)]+i\pi)&=f(z\vee h) \vee g(z\vee h) \vee ([f(z)\vee g(z)]+i\pi)  \\
   &= f(z\vee h) \vee g(z\vee h) \vee (f(z)+i\pi)\vee(g(z)+i\pi)\\
   &= f(z\vee h) \vee (f(z)+i\pi )\vee g(z\vee h)\vee (g(z)+i\pi)
\end{align*}
Next, distribute the $-h$ over the joins and take the limit as $h$
approaches $-\infty$.  The result then follows.
\item[(v)] We compute using induction and the previous result
\begin{align*}
  D(n_\vee f)&=D(f\vee \cdots \vee f)  \\
  &= n_\vee D(f).
\end{align*}
\item[(vi)] First consider the "numerator" of the "difference
quotient": {\setlength\arraycolsep{2pt}
\begin{align*}
   [f(z\vee h) + g(z\vee h)]\vee[(f(z)+g(z))+i\pi]&=&[f(z\vee h)+g(z\vee h)]\vee [f(z)+g(z)+i\pi ]  \\
   &=& [f(z\vee h)+g(z\vee h)] \vee [f(z\vee h)+g(z)] \\
   & {} & \vee [f(z\vee h)+g(z)+i\pi] \vee[f(z)+g(z)+i\pi ] \\
   &=& \left( f(z\vee h)+[g(z\vee h)\vee (g(z)+i\pi)]  \right)\\
   & & \vee \left( g(z)+[f(z\vee h)\vee (f(z)+i\pi)]   \right)
\end{align*}}
The result then follows upon distributing $-h$ and taking limits.
\item[(vii)] To show that $D(nz)=(n-1)z+\ln n$, use (vi).
\item[(viii)]To prove
$D(z^n)=\log\left(\displaystyle
\frac{nz^{n-1}e^{z^{n}}}{e^{z}}\right)$, consider the "difference
quotient" $\left[ (z\vee h)^n -h\right]\vee \left[ z^{n}+i\pi -h
\right]$, in which we have already distributed the $-h$. This
becomes
\begin{align*}
  \log \left\{ \exp \left[ (z\vee h)^{n}-h \right]+\exp \left[ z^{n}+i\pi -h \right] \right\} &= \log \left\{ \frac{\exp \left[ (z\vee h)^{n} \right]+\exp\left[ z^{n}+i\pi \right]}{\exp(h)}
  \right\}\\
  &=\log \left\{ \frac{\exp \left[ (z\vee h)^{n} \right]-\exp\left[ z^{n} \right]}{\exp(h)}
  \right\}
\end{align*}
Now we can apply l'Hospital's Rule, since the numerator and the
denominator each approach $0$ as $h\to -\infty$.  Thus, we let
$f(h)=\exp \left[ (z\vee h)^{n}-h \right]+\exp \left[ z^{n}+i\pi
-h \right]$ and $g(h)=\exp(h)$.  Then, computing ordinary
derivatives with respect to $h$ of $f$ and $g$ we obtain:
\begin{align*}
f'(h)&=\exp \left[ (z\vee h)^{n} \right]\cdot n \cdot (z\vee h)^{n-1} \cdot \frac{d}{dh}(z\vee h) \\
&=\exp \left[ (z\vee h)^{n} \right]\cdot n \cdot (z\vee h)^{n-1}
\cdot \frac{\exp(h)}{\exp(z)+\exp(h)}
\end{align*}
by Eq. ~(\ref{dervofjoin}). Further, $g'(h)=\exp(h)$. Therefore,
\begin{align*}
\frac{f'(h)}{g'(h)}&= \frac{\exp \left[ (z\vee h)^{n} \right]\cdot
n \cdot (z\vee h)^{n-1} \cdot
\displaystyle \frac{\exp(h)}{\exp(z)+\exp(h)}}{\exp(h)} \\
&=\frac{\exp \left[ (z\vee h)^{n} \right]\cdot n \cdot (z\vee
h)^{n-1}}{\exp(z)+\exp(h)}
\end{align*}
Taking the limit as $h\to -\infty$, we have $(z\vee h)^{n}\to
(z\vee -\infty)^{n})=z^{n}$ and thus,
\begin{align*}
\frac{f'(h)}{g'(h)}&=\frac{\exp \left[ (z\vee h)^{n} \right]\cdot
n \cdot (z\vee h)^{n-1}}{\exp(z)+\exp(h)} \\
& \to \frac{\exp\left[z^{n} \right]\cdot n \cdot z^{n-1}}{\exp(z)} \\
& = \frac{nz^{n-1}e^{z^{n}}}{e^{z}}.
\end{align*}
By l'Hospital's Rule $\displaystyle \frac{f(h)}{g(h)}$ approaches
the same limit as $\displaystyle \frac{f'(h)}{g'(h)}$ and, by the
continuity of the logarithm,
we obtain \\
$D(z^{n})=\log\left[ \displaystyle \frac{nz^{n-1}e^{z^{n}}}{e^{z}}
\right]=z^{n}-z + (n-1)\log(z)+\ln(n)$.
\item[(ix)] Finally, we compute
\begin{align*}
   D^{n}(nz)&=D^{n-1}D(nz)  \\
   &= D^{n-1}n_\vee (n-1)z \\
   &= n_\vee D^{n-1}(n-1)z \\
   &= n_\vee (n-1)_\vee D^{n-2}(n-2)z \\
   & ~ ~ \vdots \\
   &= n_\vee (n-1)_\vee \cdots_\vee (n-n+1)_\vee(n-n)z \\
   &= (n!)_\vee 0 \\
   &= \ln (n!) \\
   &= \ln 1 + \cdots + \ln n
\end{align*}
\end{itemize}
\end{proof}
\begin{rem}
The use of l'Hospital's Rule in the demonstration in part (viii)
of Theorem \ref{bigthm} is not restricted to that specific
instance.  For a general analytic function $f$ we consider the
following:
\begin{align*}
Df(z) &= \lim_{h \rightarrow -\infty} \left \{ \left [ f(z\vee h)
\vee (f(z)+i\pi) \right ] - h \right \} \\
&= \lim_{h \rightarrow -\infty} \left\{ \log \left[ \exp[f(z\vee h) ]+\exp(f(z)+i\pi) \right] - \log[\exp(h)] \right\} \\
&= \lim_{h \rightarrow -\infty} \log \left\{ \frac{\exp[f(z\vee h)
]+\exp[f(z)+i\pi]}{\exp(h)}\right\}
\end{align*}
Since the limit of the numerator and the denominator each approach
0, we can apply l'Hospital's Rule and take the (ordinary)
derivative of each with respect to $h$:
\begin{align*}
Df(z) &=\lim_{h \rightarrow -\infty} \log \left\{
\frac{\exp[f(z\vee h) ]+\exp[f(z)+i\pi]}{\exp(h)}\right\} \\
&= \lim_{h \rightarrow -\infty} \log \left\{
\frac{\exp[f(z\vee h)]\cdot \displaystyle f'(z\vee h)\cdot \displaystyle \frac{\exp(h)}{\exp(z)+\exp(h)}}{\exp(h)}\right\} \\
&=\lim_{h \rightarrow -\infty} \log \left\{
\frac{\exp[f(z\vee h)]\cdot \displaystyle f'(z\vee h)}{\exp(z)+\exp(h)} \right\} \\
&=\log \left\{\frac{\exp[f(z)]f'(z)}{\exp(z)} \right\} \\
&=f(z)+\log\left[f'(z)\right] - z.
\end{align*}
In other words, we have proved:
\end{rem}
\begin{thm}
If $f:G \to \mathbb{C}$ is an analytic and $\vee$-differentiable
function defined on a open subset $G\subseteq\mathbb{C}$, then for
all $z\in\mathbb{C}$,
\begin{equation}\label{mainresult}
  Df(z)=f(z)+\log[f'(z)]-z.
\end{equation}
\end{thm}
\begin{rem}
This general fact supplies some more elegant proofs of the special
formulae in Theorem \ref{bigthm}.  For example, the following is a
proof of the "product rule" for $\vee$-differentiation contained
in (iv) of that theorem:
\begin{align*}
  D(f + g)(z) &= f(z)+g(z) + \log\left[f'(z)+g'(z)\right]-z \\
  &=f(z)+g(z) + \left( \log[f'(z)]\vee \log[g'(z)]\right) -z \\
  &=\left[f(z)+g(z)+\log[f'(z)]-z\right] \vee \left[
  f(z)+g(z)+\log[g'(z)]-z\right]\\
  &= \left[Df(z) + g(z)\right] \vee \left[f(z) + Dg(z)\right].
\end{align*}
\end{rem}
\begin{rem}
Eq ~(\ref{mainresult}) also enables us to prove new results such
as
\begin{itemize}
  \item[(i)] $D[ae^{bz}]=ae^{bz}+(b-1)z+\log(ab)$
  \item [(ii)] $D[\exp(z)]=\exp(z)$
  \item [(iii)] $D[\log(z)]=-z$
\end{itemize}
It is interesting that the exponential function retains with
respect to the $\vee$-derivative the invariance it enjoys with
respect to the ordinary derivative.  Additionally, the
$\vee$-derivative of logarithm takes on the analogue value:  $-z$
is to $+$ what $1/z$ is to $\times$.

There are many more formulae analogous to classical ones that we
could prove, the chain rule being one of the more important ones.
However, space limitations prevent us from presenting them here.
\end{rem}
\section{Generalizations of the Derivative}
The approach used to define the $\vee$-derivative can be taken to
define a generalized derivative for any $n$. Let us introduce the
following new notation:  For each $n\in\mathbb{Z}$, let $\sim_n z$
denote the inverse of $z\in\mathbb{C}$ with respect to the
operation $\oplus_n$, and let $0_n$ denote the identity element
with respect to $\oplus_n$.  Thus $0_{-1}=-\infty$, $0_0=0$ and
and $0_1=1$.  Also, $\sim_{-1}z=z+i\pi$, $\sim_0 z=-z$, and
$\sim_1 z=1/z$.  Recalling that $\vee=\oplus_{-1}$, and
introducing the notation $D_{-1}$ for our new derivative, we see
that
\begin{align*}
  D_{-1}f(z)&= \lim_{h \rightarrow -\infty} \left\{ \left[ f(z\vee h) \vee (f(z)+i\pi) \right] - h \right\} \\
  &= \lim_{h \rightarrow 0_{-1}}  \left\{ \left[ f(z\oplus_{-1} h) \oplus_{-1} (\sim_{-1}f(z)) \right] \oplus_{0}
  [\sim_{0}h]
  \right\}
\end{align*}
Formally, this suggests that we attempt to generalize the
derivative as follows:
\begin{defn} For any $n\in\mathbb{Z}$, $f:G\to\mathbb{C}$ analytic on an open
subset $G\subseteq\mathbb{C}$, and $z\in G$, define the
\textit{\textbf{$n$-derivative}} of $f$ at the point $z$ as \\
\begin{equation}
D_{n}f(z):=\lim_{h \rightarrow 0_n} \left\{ \left[ f(z\oplus_{n}
h) \oplus_{n} (\sim_{n}f(z)) \right] \oplus_{n+1} [\sim_{n+1}h]
\right\}
\end{equation}
\end{defn}
\begin{rem}
By applying this to the case of $n=1$, and recalling both that
$x\oplus_2 y=\exp[\log(x)\log(y)]$ and $\sim_2
z=\exp\left[\displaystyle\frac{1}{\log(z)} \right]$, we compute:
\begin{align*}
D_1 f(z)&= \lim_{h \rightarrow 1} \exp\left\{ \log\left[
\frac{f(zh)}{f(z)} \right]  \log\left[ \exp\left(\frac{1}{\log(h)}\right)\right] \right\}\\
&=\lim_{h \rightarrow 1} \exp\left\{ \frac{\log\left[\displaystyle
\frac{f(zh)}{f(z)} \right]}{\log(h)}\right\} \\
 &=\lim_{h \rightarrow 1} \exp\left\{ \displaystyle
\frac{\log[f(zh)]-\log[f(z)]}{\log(h)}\right\} \\
\end{align*}
which is clearly reminiscent of a difference quotient.  We can
refine this further with the aid of l'Hospital's Rule as follows:
\begin{align*}
D_1 f(z)&=\lim_{h \rightarrow 1} \exp\left\{ \displaystyle
\frac{\log[f(zh)]-\log[f(z)]}{\log(h)}\right\}\\
&=\lim_{h \rightarrow 1} \exp\left[ \frac{f'(zh)zh}{f(zh)}\right]\\
&=\exp\left[ \frac{f'(z)z}{f(z)}\right].
\end{align*}
\end{rem}
 \emph{}
If this formula makes sense, then we should also be able to
generalize Eq. ~(\ref{mainresult}) in the following manner:
\begin{equation}\label{genmain}
  D_n f(z)= f(z)\oplus_{n+1}\log\left[
  D_{n+1}f(z)\right]\oplus_{n+1}[\sim_{n+1}z]
\end{equation}
This in turn should allow us to compute the ordinary
$0$-derivative $D_0$ in terms of the next higher $1$-derivative
$D_1$.  In other words,
\begin{align*}
f'(z)=D_0f(z)&=f(z)\oplus_1 \log\left[
  D_{1}f(z)\right]\oplus_{1}[\sim_{1}z]\\
  &=f(z)\cdot \log\left[
  D_{1}f(z)\right]\cdot\frac{1}{z}\\
  &=f(z)\cdot \log\exp\left[
  \frac{f'(z)z}{f(z)}\right]\cdot\frac{1}{z}\\
  &=f(z)\cdot\frac{f'(z)z}{f(z)}\cdot\frac{1}{z}\\
  &=f'(z).
\end{align*}

We note that the exponential function is again invariant with
respect to $D_1$, just as it is with respect to $D_{-1}$ and
$D_0=\displaystyle \frac{d}{dz}$. For,
\begin{align*}
D_{1}\exp(z)&=\exp\left[\frac{\displaystyle \frac{d}{dz}\exp(z)\cdot z}{\exp(z)}\right]\\
&=\exp\left[\frac{\exp(z)\cdot z}{\exp(z)}\right]\\
&=\exp\left[z \right]
\end{align*}

We are tempted to conjecture that the exponential function is
probably invariant with respect to every $D_n$ as defined in Eq.
~(\ref{genmain}).  However, we choose to leave this for subsequent
publications.

 Naturally, these results are primarily formal at this
point and need to be put on more solid analytical footing by
careful consideration of domains of definition and branches of the
logarithm.  However, their formal consistency is very encouraging
and seems to indicate that the generalized derivative, along with
the generalized operations, yields a countably infinite hierarchy
of structures on which to conduct complex analysis.

\section{Conclusions and Prospects for Further Research}
In this paper I have presented a natural chain of binary
operations on the reals and complex numbers.  The ordinary field
operations of addition and multiplication are seen to be members
of that chain. Indeed, the chain is constructed from the operation
of addition using combinations of exponential and logarithmic
mappings.   In particular, we have explored one of these new
operations -- the join -- and have seen that a new kind of
difference quotient and derivative can be defined using it.  The
formal methods used in defining the derivative appear to be easily
generalizable to any pair of adjacent operations in our chain. The
main result is the establishment of a formula that relates in
essence any one of these derivatives to that of the ordinary
calculus.  Indeed, the ordinary derivative is seen to be a special
case of the more general derivative.

It should be apparent that we have in this theory the beginnings
of a potentially fruitful field of study.  A number of the results
are sufficiently compact and appealing to warrant further
exploration in this area.

It is important to realize that the each operation in the chain is
fully compatible with the existing topological and algebraic
structure of the field of complex numbers.  This is due to the
continuous and smooth nature of the logarithm and exponential
functions.  Because these new operations can be related directly
to the classical operations of addition and multiplication, they
present not merely an alternative representation of the legacy
field structure but rather an enrichment of it.  It is this key
fact that seems to hold out so much promise for further results
and new insights into old problems.  It also seems to indicate
that we can break out of the two-operation mold of traditional
ring theory and begin to explore more intricate structures with
many more intimately related operations.

In subsequent investigations I will focus on exploring in more
depth the relationship between the new derivative and the
classical one and on defining integration.  In particular, I will
show that the generalized derivative as given by Eq.
~(\ref{genmain}) is indeed a ring derivation for each $n$. I will
also attempt to gain more insight into the geometric
interpretation of the join and the $\vee$-derivative in the spirit
of Needham\cite{Nee97}. Furthermore, I will attempt to develop a
systematic theory of infinite joins analogous to that of infinite
series and to exhibit some new infinite join representations of
various important functions and constants.

The relationship between these natural operations and the $\max$
operation and max-plus rings and algebras in general would also
seem to be of interest. It is easy to prove, for example, that
$\max\{x,y\}\leq x\oplus_n y$ for every $n\in\mathbb{Z}$.  It
appears likely that $\max\{x,y\}= \displaystyle \lim_{n
\rightarrow -\infty} x\oplus_n y$.  However, whereas $\max$ is not
a smooth operation, $\oplus_n$ is. This might be of benefit in
applications to nonlinear systems theory where differentiability
is always a desirable quality and the potential for transforming
non-linear problems into linear ones is always attractive.

The territory opened up by this research topic appears to be
potentially vast and certainly much larger than can be adequately
explored by a single researcher.  It is my hope that this paper
will spark the interest of other researchers and that
the number of new results and applications of this field will grow
rapidly.

% ------------------------------------------------------------------------

\subsection*{Acknowledgment}
I wish to express my gratitude to my wife Leandra for her
encouragement and patience while I pursued this study and to Randy
Jaffe, my supervisor at BEI Systron Donner Inertial Division, who
encouraged and kindly reviewed this work.
% ------------------------------------------------------------------------
%GATHER{Xbib.bib}   % For Gather Purpose Only
%GATHER{Paper.bbl}  % For Gather Purpose Only
\bibliographystyle{amsplain}
\bibliography{xbib}
\end{document}